\documentstyle[12pt]{article}

\def\ddd#1{}  
\def\xxx#1{}	
\def\mat{\def}	\def\grk{\def}	\def\pol{\def}
	
\xxx{
\setlength{\textwidth}{164mm}	\setlength{\oddsidemargin}{2mm}
\setlength{\textheight}{230mm}	\addtolength{\topmargin}{-2cm}
}

\mat\0#1{\makebox{\(#1\)}} 
\mat\8#1{\makebox{\(\displaystyle #1\)}} 

\mat\9#1{\0{ \left( {#1} \right) }}
\mat\1#1{\0{ \left| {#1} \right|  }}	
\mat\2#1{\0{ \left\|{#1} \right\| }}
\mat\3#1{\0{ \left\langle {#1} \right\rangle}}	
\mat\4#1{\0{ \left[ {#1} \right]  }}	
\mat\5#1{\0{ \left\{{#1} \right\} }}	
\mat\6#1#2{\0{#1\in #2}}
\mat\7#1#2#3{\0{#1:#2\to #3}}	

\pol\cu{\raisebox{-1.45ex}[-2ex][-2ex]{``}}
\pol\polae#1#2#3{\makebox[0em][l]{\hspace*{#1}\raisebox{#2}[0mm][0mm]{`}}#3}
\pol\a.{\polae{.31em}{-1.45ex}{a}}
\pol\A.{\polae{.44em}{-1.45ex}{A}}
\pol\e.{\polae{.20em}{-1.45ex}{e}}
\pol\E.{\polae{.41em}{-1.45ex}{E}}
\pol\u.{{\l}}

\mat\AA{{\cal A}}	\mat\JJ{{\cal J}}	\mat\SS{{\cal S}}	
\mat\BB{{\cal B}}	\mat\KK{{\cal K}}	\mat\TT{{\cal T}}	
\mat\CC{{\cal C}}	\mat\LL{{\cal L}}	\mat\UU{{\cal U}}	
\mat\DD{{\cal D}}	\mat\MM{{\cal M}}	\mat\VV{{\cal V}}	
\mat\EE{{\cal E}}	\mat\NN{{\cal N}}	\mat\WW{{\cal W}}	
\mat\FF{{\cal F}}	\mat\OO{{\cal O}}	\mat\XX{{\cal X}}	
\mat\GG{{\cal G}}	\mat\PP{{\cal P}}	\mat\YY{{\cal Y}}	
\mat\HH{{\cal H}}	\mat\QQ{{\cal Q}}	\mat\ZZ{{\cal Z}}	
\mat\II{{\cal I}}	\mat\RR{{\cal R}} 

\mat\N{{\bf N}}		\mat\Z{{\bf Z}}		
\mat\R{{\bf R}}		\mat\C{{\bf C}}

\grk\Gam{\Gamma}	\grk\Del{\Delta}	\grk\Th{\Theta}   
\grk\Lam{\Lambda}       \grk\Sig{\Sigma}	\grk\Ups{\Upsilon}	
\grk\Yps{\Upsilon}    	\grk\Ph{Phi}		\grk\Ps{\Psi} 
\grk\Om{\Omega}     

\grk\al{\alpha}	 	\grk\gam{\gamma}	\grk\del{\delta}	
\grk\eps{\varepsilon} 	\grk\th{\vartheta}	\grk\kap{\kappa}	
\grk\lam{\lambda}	\grk\rh{\varrho}	\grk\sig{\sigma}	
\grk\ups{\upsilon} 	\grk\ph{\varphi}	\grk\ps{\psi} 		
\grk\om{\omega}

\xxx{
\grk\alfa{\alpha} 	\grk\dzeta{\zeta} 	\grk\dz{\zeta} 
\grk\teta{\vartheta}	\grk\jota{\iota} 	\grk\mi{\mu}	
\xxx{\grk\ni{\nu}} 	
\grk\ksi{\xi} 		\grk\ro{\varrho} 	
\grk\yps{\ypsilon}
} 

\newtheorem{xDef}{Definition}
\newtheorem{xTheo}{Theorem}
\newtheorem{xCor}[xTheo]{Corollary}
\newtheorem{xLem}[xTheo]{Lemma}
\newtheorem{xProp}[xTheo]{Proposition}
\newtheorem{xEx}{Example}
\newtheorem{xQue}[xTheo]{Question}

\newenvironment{Definition}{\begin{xDef}\rm }{\end{xDef}}
\newenvironment{Theorem}{\begin{xTheo}\sl }{\end{xTheo}}
\newenvironment{Corollary}{\begin{xCor}\sl }{\end{xCor}}
\newenvironment{Lemma}{\begin{xLem}\sl }{\end{xLem}}

\mat\Pf{{\em Proof. }}       	
\mat\fP{\rule{2.5mm}{2.5mm}\vspace*{1.5ex}}
\mat\barr{\begin{array}} 	\mat\earr{\end{array}} 
\mat\bequ{\begin{equation}} 	\mat\eequ{\end{equation}} 
\mat\bequarr{\begin{eqnarray}}  \mat\eequarr{\end{eqnarray}}
\mat\bite{\begin{itemize}} 	\mat\eite{\end{itemize}}

\mat\Then{\Longrightarrow}
\mat\to{\rightarrow}		\mat\from{\leftarrow}
\mat\To{\longrightarrow}	\mat\From{\longleftarrow}
\mat\onto{\stackrel{{\scriptsize\rm onto}}{\to}}

\mat\stk{\stackrel}		\mat\fr#1#2{{\textstyle\frac{#1}{#2}}}
\mat\xl#1{\mbox{\rm\ \ #1\ \ }}
\mat\r#1{\mbox{\rm #1}}
\mat\til#1{\0{\widetilde{#1}}}	\mat\ovl#1{\0{\overline{#1}}}

\mat\d.{\partial}	
\mat\c.{\cdot}		\mat\x.{\times}	

\mat\eq{\equiv}		\mat\app{\approx} 
\mat\ctin{\subset}      \mat\cts{\supset}
\mat\mns{\setminus} 	
\mat\Union{\Bigcup}     \mat\Inters{\Bigcap}

\def\xto{ 
\setlength{\unitlength}{1ex} 
	\begin{picture}(5,1.6)
\put(0.8,0.8){\vector(1,0){4}}
	\end{picture}
}

\def\xline{ 
\setlength{\unitlength}{1ex} 
	\begin{picture}(5,1.6)
\put(0.8,0.8){\line(1,0){4}}
	\end{picture}
}

\begin{document} \sloppy

\title{The three-colour hat guessing game on the cycle graphs}

\author{Witold W.\ Szczechla 
\thanks{Department of Mathematics, Informatics and Mechanics, 
University of Warsaw, ul.Banacha\,2, 02-097 Warszawa, Poland.
{\tt witold@mimuw.edu.pl}
} 
} 

\date{} 

\maketitle 

\begin{abstract} We study a cooperative game in which each member of a 
team of $N$ players, wearing coloured hats and situated at the vertices of 
the cycle graph $C_N$, is guessing their own hat colour merely on the 
basis of observing the hats worn by their two neighbours without 
exchanging the information.  Each hat can have one of three colours. A 
predetermined guessing strategy is winning if it guarantees at least one 
correct individual guess for every assignment of colours.  We prove that a 
winning strategy exists if and only if $N$ is divisible by $3$ or $N=4$.  
This problem represents an example of a relational system using incomplete 
information about an unpredictable situation, where at least one 
participant has to act properly. \end{abstract}
 

\section{Introduction} 

$N$ ladies wearing white hats are sitting around the table and discussing 
a tricky task which is going to be presented to them by the Wizard.  They 
know he will suddenly paint each hat one of three colours (green, orange 
or purple) in an unpredictable way and then ask each of them to 
independently guess her own hat colour.  The light is so dim that everyone 
will only see the hat colours of her two neighbours.  If at least one of 
the ladies guesses right, they will all win; if they all guess wrong, they 
will lose; and they want to be absolutely certain of winning. However, 
can they devise a winning strategy before inviting the Wizard? The answer, 
depending on the number $N$, is presented in this paper. 

\subsection{Motivation} 

This game is an example of a relational system which uses incomplete 
information about an unpredictable situation, where at least one link has 
to act or work properly. 
Problems of this kind have become popular in recent years both as 
mathematical puzzlers (see \cite{NYTimes}, \cite{Winkler}) and research 
subjects.  Basic results so far concerned two colours (instead of three), 
or unrestricted visibility (corresponding to a complete graph), or 
probabilistic variants (the expected number of the correct guesses): see 
\cite{Ebert}, \cite{Feige1}, \cite{Feige2}, \cite{Riis}, \cite{WCR}, 
\cite{Krzyw0}, \cite{Krzyw1} and overviews in \cite{BHKL} and 
\cite{Krzyw2}. The round-table problem described above has until now 
remained open for all $N>5$.


\subsection{Formalism} 
  
The team players are seeing each other along the edges of the cyclic graph 
$C_N$.  Let the set $V_k=\5{v_1(k),v_2(k),v_3(k)}$ represent the three 
different appearances of the $k$-th hat, where $k$ is counted 
modulo $N$ in  the positive direction (to the right). It will be 
technically convenient to regard them either as pairwise disjoint $(0\leq 
k<N)$, or simply as $V_k=\Z_3$. A cyclic notation will also be used:   
\[ 
	\begin{array}{ll}
 &  v_i^k=v_j^m=v_i(k)=v_j(m)~~\r{and}~~V_k=V_m, \\  
\r{where}  &  
i,j,k,m\in\Z, ~~ i\equiv j\,(\r{mod}\,3), ~~ k\equiv m\,(\r{mod}\,N). 
	\end {array}
\] 
An individual guessing strategy of Player\,$k$ is represented by a 
function
	\[
	\7{f_k}{V_{k-1}\times V_{k+1}}{V_k}, 
	\]
which may simply be regarded as a function 
$(i,j)\mapsto r$, where $f_k(v_i^{k-1},v_j^{k+1})=v_r^k$.
A {\em composite strategy} is a sequence 
	\[ f=\9{f_1,\ldots,f_N}, 
	\] 
or equivalently, a function $\Z\ni k\mapsto f_k$, satisfying 
$f_{k+N}=f_k$.  

According to the assumed rules, strategy $f$ is {\em winning\/} if and 
only if there is no sequence \9{s_1,s_2,\ldots,s_N} satisfying
	\[  
\6{s_k}{V_k} \r{~~and~~} s_k\neq f_k(s_{k-1},s_{k+1}) 
\r{~~for~~}  k=1,\ldots, N, \r{~~where~~} s_0=s_N,\ s_{N+1}=s_1.  
	\] 
If there exists such a sequence, $f$ is a {\em losing\/} strategy. 


\subsection{Hat games on graphs} 

In a more general setting an arbitrary `visibility' pattern can be 
assumed.  
For a broader exposition, see \cite{BHKL} and \cite{Krzyw2}.  
The directed `visibility graph' $\Gam$ has $N$ vertices 
corresponding to the players, and edges $\vec{AB}\in E(\Gam)=E$ wherever   
player~$A$ is seen by player~$B$.  For each vertex $v\in V(\Gam)=V$ a 
nonempty set of `colours' $V_v$ is known to all.  For each 
`assignment of colours', i.e., a selector \7{g}{V}{\bigcup_v V_v} 
with $g(v)\in V_v$, each player $u\in V$ tries to guess $g(u)$ by 
using a function 
\[  
\7{f_u}{\prod V_v}{V_u} ~~~
(\r{product taken over} ~~ \vec{vu}\in E).  
\] 
as an individual strategy.  The combined, or collective, strategy is the 
collection $f=\5{f_u:\,u\in V}$. The game is thus played against an opponent 
assigning the colours (the Wizard, the Demon, Chance,  
etc.,~in a fantasy world). In this paper, the notion of winning or losing 
refers to the cooperative players. The strategy effectiveness depends only 
on the numbers of possible colours, i.e., the function $h$ given by 
$h(v)=|V_v|$. Let $X_h(f,g)$ denote the number of correct guesses.  
The deterministic minimax approach defines the {\em value} of this game as  
	\[
  \mu(h) = \mu(\Gam,h) = \max_f\,\min_g X_h(f,g).  
	\] 
This paper concerns the minimal condition $\mu(h)>0$ where   
$f$ is a winning strategy if 
	\[
  \min_g X_h(f,g)>0.
	\] 
J.Grytczuk has conjectured that a winning strategy exists provided 
$|V_v|\leq \r{deg}_-(v)$ for every \6{v}{V(\Gam)}.  One consequence of our 
main result is that the weaker condition $|V_v|\leq \r{deg}_-(v)+1$ is 
generally not sufficient.


\subsection{Examples with $N<5$}\label{Sec34} 

\paragraph{The game on \mbox{\boldmath $C_2$}} 
In the simplest puzzle (outside our main problem, though) there are just 
two players and two possible hat colours.  In this situation one person 
should guess that their hats have the same colour and the other person 
sholud guess the opposite.  If one interpretes the colours as elements 
$A,B\in\Z_2$, the effect can be written as an alternative: 
	\[
A=B\ \ \r{or}\ \ B=A+1. 
	\] 

Next, suppose Player\,1 hat can still have two colours, but Player\,2 hat 
can have three colours.  Then there are six possible colour assignments.  
With any strategy, Player\,1 guesses right for three assignments, 
Player\,2 for two.  Since the number of assignments is $6>3+2$, they can 
both be wrong and they have no winning strategy. (In the above cases, the 
players might as well guess the other person's colour while knowing their 
own.) 


\paragraph{An algebraic strategy for \mbox{\boldmath $C_3$}} 
If  $A,B,C\in\Z_3$ represent the appearances of hats, then a 
winning strategy can be based, for instance, on the alternative:
	\[
A=-B-C\ \ \r{or}\ \ B=-C-A-1\ \ \r{or}\ \ C=-A-B+1, 
	\]
clearly valid in $\Z_3$. 


\paragraph{An algebraic strategy for \mbox{\boldmath $C_4$}} Let variables 
$A,B,C,D\in\Z_3$ represent elements of the sets $V_1,V_2,V_3,V_4$, 
respectively.  Then a winning strategy $f$ can be based on the 
following alternative:
	\begin{equation}\label{eq_4} 
	\left\{ \begin{array}{cccc} 
	  & A & = &  D+B  \\ 
\r{or}    & B & = & -A-C  \\ 
\r{or}    & C & = &  B-D  \\ 
\r{or}    & D & = &  C-A. 
	\end{array} \right. 
	\end{equation} 
To verify (\ref{eq_4}), let us suppose the first and third equalities are 
false.  In $\Z_3$ this implies
	\begin{equation} \nonumber  
	\left\{ 
	\begin{array}{ccc} 
D+B  & = &  A \pm 1  \nonumber \\ 
B-D  & = &  C \pm 1. \nonumber 
	\end{array} 
	\right. 
	\end{equation} 
If the signs above are opposite, we add the equations to get $2B=A+C$, 
equivalent to the second equation of (\ref{eq_4}).  If the signs are the 
same, we subtract the equations to get $2D=A-C$, equivalent to the last 
equation of (\ref{eq_4}).  Thus indeed, strategy $f$ wins.


\section{The main result} 

Let us observe the following (inconvenient, as it turns out) property of 
the last two examples.  For any $b\in V_m$ and $c\in V_{m+1}$ 
there is exactly one $d\in V_{m+2}$ satisfying $f_{m+1}(b,d)=c$ and 
exactly one $a\in V_{m-1}$ satisfying $f_m(a,c)=b$.  However, winning 
strategies with this property are {\em not\/} possible for $N>4$ (see 
Section~\ref{Sec0}).  Our method will thus produce another kind of 
strategy for some $N$-gons whenever possible, while demonstrating the 
losing case for all the rest. The main result of this paper is:

\begin{Theorem}\label{T1} 
In the three-colour hat guessing game on the cycle 
of length $N$ a winning strategy exists if and only if $N$ is divisible 
by three or $N=4$.  
\end{Theorem}


\subsection*{Proof of the main result} 

An interplay of various relatively simple local and global combinatorial 
methods will be used. 

\subsection{Admissible paths in the enlarged graph}  

Let us introduce a larger graph $G=G_N$ 
(which we will also denote $3\ast C_N$) 
whose $3N$-element set of vertices is 
$V=V(G)=\bigcup_{k=1}^N V_k$, 
and $9N$-element set of edges is 
	\[ E=E(G)= 
\{ \ovl{v_i(k-1)v_j(k)}:\ k=1,\ldots, N;\ i,j=1,2,3 \}. 
	\]

\paragraph*{Remark} 
An analogous construction can be applied to any visibility 
graph $\Gam$ and any {\em height function\/} \7{h}{V(\Gam)}{\N\mns\5{0}} 
(whose values are the numbers of possible colours).  The resulting graph, 
which may be denoted $G=\ast\Gam$, has
	\[
V(G) = \5{ (i,v):\, v\in V(\Gam),\, i=1,\ldots,h(v)}
	\]
and 
	\[
E(G) = \5{ \vec{(i,v)(j,u)}:\, \vec{vu}\in E(\Gam)}. 
	\]

Now let us consider a (composite) strategy $f$. 
	\begin{Definition} 
Let $J$ be any set of consecutive integers.  A path $(s_k)_{k\in J}$ in 
the graph $G$ will be called {\em $f$-admissible\/} (or 
simply {admissible\/}, when $f$ is fixed) if
	\[ s_k\in V_k\ \ \r{for}\ \ k\in J 
	\]  
and 
	\[ s_k\neq f_k(s_{k-1},s_{k+1})\ \ \r{whenever}\ \ k-1,k+1\in J.  
	\]  
$\Box$
	\end{Definition} 

Thus, a path is admissible if and only if all its $2$-edge segments (i.e., 
sub-paths of length~$2$) are admissible.  It is clear that strategy $f$ is 
winning if and only if the graph $G$ contains no $f$-admissible path (of 
infinite length) which is periodic and has period $N$, or equivalently, no 
$f$-admissible path of length $N+1$ whose last edge coincides with the 
first.  However, the definition does not direcly settle the 
question of whether any periodic admissible path exists, and if it does, 
whether it can have period $N$ (or at least less than $9N$).


The set of all the $f$-admissible paths (of all lengths) will be denoted 
$\AA(f)$.  The set of all the edges between $V_k$ and $V_{k+1}$ will be 
denoted by 
	\[ E_{k,k+1} = V_k\times V_{k+1} 
	= \5{\ovl{bc}:\,\6{b}{V_k}~\r{and}~\6{c}{V_{k+1}}}. 
	\]   
	\begin{Definition} Let $f$ be a fixed strategy. 
For any edge \6{\ovl{bc}}{E_{k,k+1}}, define 
\[  \ell_+(\ovl{bc}) = \#\5{d\in V_{k+2}:~\ovl{bcd}\in\AA(f)}, \] 
i.e., the number of the immediate admissible continuations of $\ovl{bc}$ 
to the right, and 
\[  \ell_-(\ovl{bc}) = \#\5{a\in V_{k-1}:~\ovl{abc}\in\AA(f)}, \] 
i.e., the number of the analogous continuations to the left.  $\Box$
	\end{Definition}
Hence, in general, we have $\ell_+(\ovl{bc}), \ell_-(\ovl{bc}) \in 
\5{0,1,2,3}$.


\begin{Lemma} \label{L1} Consider a fixed strategy $f$. 
	\begin{itemize} 
	\item[{(a)}] 
The average value of $\ell_-$ (resp.~$\ell_+$) over any three 
right-adjacent (resp.~left-adjacent) edges of $G$ equals~$2$.  
That is, for any vertex \6{b}{V_k} we have
	\[ 
\sum_{a\in V_{k-1}}\ell_-(\ovl{ab}) =  
\sum_{c\in V_{k+1}}\ell_+(\ovl{bc}) = 6
	\] 
	\item[{(b)}] If two edges of $G$ have the same left 
(resp.~right) endpoint then one of them has at least two admissible 
immediate continuations to the right (resp.~left).  That is, for any 
vertex \6{b}{V_k} and any two distinct vertices \6{c_i}{V_{k+1}} \9{i=1,2} 
there is a choice of\/ \6{i}{\5{1,2}} and two distinct vertices 
\6{d_1,d_2}{V_{k+2}} such that $\ovl{bc_id_j}\in\AA(f)$ for $j=1,2$; 
the analogous fact holds for passages to the left.
	\item[{(c)}] If the graph $G$ contains an $f$-admissible path 
$\ovl{s_1\ldots s_n}$ such that $2\leq n\leq N-1$ and 
	\[ \ell_-(\ovl{s_1s_2})+\ell_+(\ovl{s_{n-1}s_n})\geq 5, \] 
then $f$ is a losing strategy.  
	\item[{(d)}] If $f$ is a winning strategy, then for every edge 
$\beta\in E(G)$ we have 
	\[ \ell_+(\beta)+\ell_-(\beta)=4.  \]
	\end{itemize} 
\end{Lemma}


\Pf (a): Consider $\ell_+$. For any vertex $d\in V_{k+2}$, the set 
$V_{k+1}$ contains two vertices different from $f_{k+1}(b,d)$, defining 
two admissible connections of~$b$ with each of the three choices of~$d$. 
(The situation with $\ell_-$ is symmetric.)

(b): For any $d\in V_{k+2}$ we can choose an \6{i}{\5{1,2}} such that 
$f_{k+1}(b,d)\neq c_i$.  Since $d$ takes three values, two of them must 
correspond to the same choice of $i$.

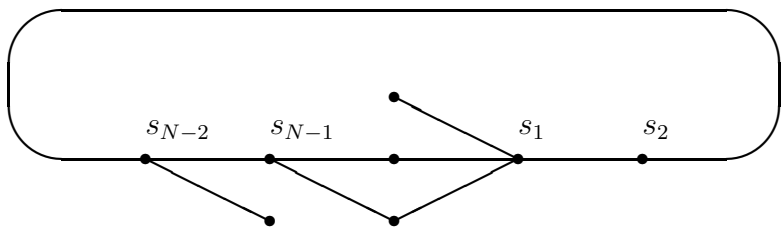
\begin{figure} \label{f23}
	\begin{center} \setlength{\unitlength}{0.4em}
	\begin{picture}(64,24)(-32,-6)   \thicklines 
   \put(-25,0){\line(1,0){50}}  
\multiput(-20, 0)(10,0){5}{\circle*{0.9}}  	
\multiput(-10,-5)(10,0){2}{\circle*{0.9}}
		\put(0,5){\circle*{0.9}}   
     \multiput(-20,0)(10,0){2}{\line(2,-1){10}}  
     \put(10,0){\line(-2, 1){10}}
     \put(10,0){\line(-2,-1){10}}  
\put(-20,2){$s_{N-2}$} 	\put(-10,2){$s_{N-1}$}
\put( 10,2){$s_1$}	\put( 20,2){$s_2$} 
     \put(-25,12){\line(1,0){50}} 		
     \put(-25,6){\oval(12,12)[l]}  
     \put( 25,6){\oval(12,12)[r]}
	\end{picture} 
\caption{Closing the path of Lemma\,2\,(c).}  
	\end{center}
	\end{figure}

(c): We may assume $\ell_-(\ovl{s_1s_2})=3$ and 
$\ell_+(\ovl{s_{n-1}s_n})\geq 2$, with $s_1\in V_1$.  By~(b), the path can 
be continued to the right until $n=N-1$.  Then the paths of the form 
$\ovl{xs_1s_2\ldots s_{N-1}y}$ are in $\AA(f)$ for all three values 
of $x\in V_0$ and at least two values of $y\in V_N=V_0$.  Now it is enough 
to choose $y\neq f_0(s_{N-1}s_1)$ to make the ends meet, obtaining an 
$N$-periodic $f$-admissible path $\ovl{ys_1s_2\ldots s_{N-1}ys_1\ldots}$.  
(This argument is partly illustrated in Figure~\ref{f23}.) 

(d): Denote $\ell(\gamma)=\ell_+(\gamma)+\ell_-(\gamma)$ for all 
$\gamma\in E(G)$.  If $\ell(\beta)>4$ for some $\beta\in E(G)$, then 
$f$ is losing by (c) applied to the single edge $\beta$.  However, (a) 
implies that the average value of $\ell(\gamma)$ over $\gamma\in 
E_{k,k+1}$ equals~$4$.  Hence, if there was an edge $\alpha\in E_{k,k+1}$ 
with $\ell(\alpha)<4$, there would also be an edge $\beta\in E_{k,k+1}$ 
with $\ell(\beta)>4$, the case already excluded. \fP

\subsection{The three categories of edges} 
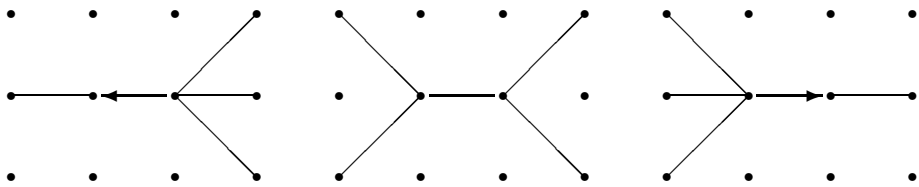
\begin{figure} \label{f3cat}
	\begin{center} \setlength{\unitlength}{0.6ex} 
	\begin{picture}(112,22)(-1,-1)  
\multiput(0,20)(10,0){12}{\circle*{1}}
\multiput(0,10)(10,0){12}{\circle*{1}}
\multiput(0, 0)(10,0){12}{\circle*{1}}
\thicklines 	\put( 19,10){\vector(-1,0){8}}
	 	\put( 51,10){\line( 1,0){8}}
	 	\put( 91,10){\vector( 1,0){8}}  
\thinlines  	
\multiput(20,10)(40,0){2}{\line(1,1){10}}	
\multiput(40,20)(40,0){2}{\line(1,-1){10}}
\multiput( 0,10)(20,0){2}{\line(1,0){10}}  
\multiput(80,10)(20,0){2}{\line(1,0){10}}  
\multiput(20,10)(40,0){2}{\line(1,-1){10}}	
\multiput(40, 0)(40,0){2}{\line(1,1){10}}
	\end{picture} 
\caption{Examples of edges (red, blue, yellow) with their admissible 
continuations.} 
	\end{center}
	\end{figure}

Let us assume that strategy $f$ 
satisfies 
	\begin{equation}\label{eq_c}
\ell_+(\gamma)+\ell_-(\gamma)=4 ~~\r{for all}~~ \6{\gamma}{E(G)}. 
	\end{equation} Then all the edges $\gamma\in E(G)$ can be divided 
into three categories:
	\begin{itemize}
	\item  If $\ell_-(\gamma)=3$ and $\ell_+(\gamma)=1$, let 
us paint $\gamma$ \underline{\em yellow\/} 
and direct it \underline{\em right.}  
	\item  If $\ell_-(\gamma)=1$ and $\ell_+(\gamma)=3$, let 
us paint $\gamma$ \underline{\em red\/} 
and direct it \underline{\em left.} 
	\item  If $\ell_-(\gamma)=\ell_+(\gamma)=2$, let us paint 
$\gamma$ \underline{\em blue\/} 
and leave it \underline{\em undirected.} 
	\end{itemize} 
The three patterns can thus be shown as in Figure~\ref{f3cat}.

\begin{Definition} Any strategy $f$ satisfying (\ref{eq_c}) will be called 
{\em balanced\/} or {\em colourable\/}. $\Box$
\end{Definition}

By Lemma~\ref{L1}(d), every winning strategy is colourable.  However, not 
all balanced strategies will be winning. 
The sets of all the yellow, red, and blue edges in 
$E(G)$ (or in $E_{k,k+1}$) will be denoted $E^+$, $E^-$, and $E^0$ (or 
$E^+_{k,k+1}$, $E^-_{k,k+1}$, and $E^0_{k,k+1}$), respectively.


\begin{Lemma}\label{L2}  If strategy $f$ is colourable, then: 
	\begin{itemize} 
\item[{(a)}] For each $k$, there are equal numbers of yellow and 
red edges in the set $E_{k,k+1}$ (i.e., $|E^+_{k,k+1}|=|E^-_{k,k+1}|$). 
	\item[{(b)}] Any three edges of $G$ having a common left   
or right end-point (i.e., left- or right-adjacent) either have three 
different colours or all are blue. 
	\item[{(c)}] If $f$ is a winning strategy and $N\geq 4$, then 
every directed edge is admissibly continued in its direction by an edge of 
the same direction.  That is, if \6{\beta}{E_{k,k+1}} is yellow, 
\6{\gamma}{E_{k+1,k+2}} and the path $\ovl{\beta\gamma}$ is 
$f$-admissible, then $\gamma$ is also yellow.  Analogously, if 
\6{\beta}{E_{k,k+1}} is red, \6{\alpha}{E_{k-1,k}} and 
$\ovl{\alpha\beta}\in\AA(f)$, then $\alpha$ is also red.
	\item[{(d)}] If $f$ is a winning strategy and $N\geq 4$, then 
every directed edge is a continuation of three edges of three different 
colours.  That is, if $\beta=\ovl{bc}\in E_{k,k+1}$ is yellow, then among 
the three edges $\ovl{ab}$ with \6{a}{V_{k-1}} one is yellow, one is red, 
and one is blue.  Analogously, if $\beta=\ovl{bc}$ is red, then among the 
edges $\ovl{cd}\in E_{k+2}$ one is in $E^+$, another in $E^-$, and the 
third in $E^0$.
	\item[{(e)}] If $f$ is a winning strategy and $N\geq 4$, then any 
periodic $f$-admissible path has one colour. Conversely (under the same 
assumption), any path of a fixed direction (red or yellow) is admissible, 
and an undirected path (blue) is admissible provided that all its vertices  
are incident to some directed edges. 
	\end{itemize}
\end{Lemma}


\Pf (a): By Lemma~\ref{L1}(a), we have 
	\[
\sum_{\gamma\in E_{k,k+1}} \ell_-(\gamma) = 
\sum_{\gamma\in E_{k,k+1}} \ell_+(\gamma) = 18. 
	\] 
The terms equal to $1$ and $3$ in the first sum correspond 
to the terms equal to $3$ and $1$, respectively, in the second.

(b): This follows from Lemma~\ref{L1}(a)(d), since the number 
$6$ can be expressed as an unordered sum of three terms equal to $1$, 
$2$ or $3$ in just two ways: $1+2+3$ and $2+2+2$. 

(c): If $\beta$ is yellow (i.e., $\ell_-(\beta)=3$) and $\gamma$ is not, 
then $\ell_+(\gamma)\geq 2$.  Then Lemma~\ref{L1}(c) applied to the path 
$\ovl{\beta\gamma}$ (where $n=3\leq N-1$) implies that $f$ is a losing 
strategy, contrary to our assumption. 

(d): Let $\beta=\ovl{bc}\in E_{k,k+1}$ be yellow.  By (b), some edge 
$\gamma=\ovl{bc'}\in E_{k,k+1}$ must be red.  Then, by (c), the left 
continuation of $\gamma$ into $E_{k-1,k}$ is also red, showing that not 
all edges $\ovl{ab}\in E_{k-1,k}$ are blue.  By (b), these edges must be 
of three colours.  

(e): Consider any path containing a yellow (resp.\ red) edge~$\beta$. By 
(c) and the definition of colouring, the edge $\beta$ has a unique forward 
(resp. backward) admissible continuation, consisting of edges of the same 
direction.  This proves that every periodic admissible path must have one 
direction or be undirected.  Conversely, any directed path is admissible 
by~(c).  

Finally, consider a blue path $\ovl{abc}$ (of length~$2$) with vertex $b$ 
incident to some directed edge.  We may suppose an edge $\ovl{bc'}$ is 
directed.  By (b) and (c), there is some left-directed edge $\ovl{bc''}$ 
uniquely continued by another edge $\ovl{a'b}\in E^-\not\ni\ovl{ab}$.  
Since $\ovl{abc''}$ is not $f$-admissible and $\ovl{ab}$ is blue, both 
remaining right continuations of $\ovl{ab}$ must be $f$-admissible, 
including $\ovl{abc}$.  \fP


	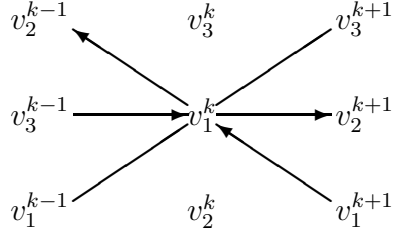
\begin{figure} \label{fX}
	\begin{center} \setlength{\unitlength}{0.7ex} 
	\begin{picture}(32,22)(-16,-11)  \thicklines 
\put(-17,9){\makebox(0,0)[b]{$v_2^{k-1}$}}  
    \put(0,9){\makebox(0,0)[b]{$v_3^k$}}  
	\put(17,9){\makebox(0,0)[b]{$v_3^{k+1}$}} 
\put(-17,-1){\makebox(0,0)[b]{$v_3^{k-1}$}} 
    \put(0,-1){\makebox(0,0)[b]{$v_1^k$}}  
	\put(17,-1){\makebox(0,0)[b]{$v_2^{k+1}$}} 
\put(-17,-11){\makebox(0,0)[b]{$v_1^{k-1}$}} 
    \put(0,-11){\makebox(0,0)[b]{$v_2^k$}}  
	\put(17,-11){\makebox(0,0)[b]{$v_1^{k+1}$}}  
\multiput( 1.5,0)(-15,  0){2}{\vector(1,0){12}}
\multiput(-1.5,1)( 15,-10){2}{\vector(-3,2){12}}
\multiput( 1.5,1)(-15,-10){2}{\line(3,2){12}}
	\end{picture} 
\caption{An example of the typical configuration at the head or tail of 
any directed edge.} 
	\end{center}
	\end{figure}

\paragraph{Alternative arguments} (b) implies (a), since there must be equal 
numbers ($0$ or $1$) of red and yellow edges left-incident to every vertex of 
$G$. Another way of proving (d) is using (c) to continue $\beta$ to the right 
with yellow edges until the last one points to the beginning of another, which 
must be $\beta$ (the first one), as two yellow edges cannot be 
(right-)incident, by (b).


\subsection{The characteristic number of a winning strategy} 

	\begin{Corollary}\label{C1} Let $f$ be a winning strategy and 
$N\geq 4$.  
	\item[{(a)}] Every directed edge $\beta\in E_{k,k+1}$ meets 
exactly two edges of $G$ having the same direction.  Moreover, one of them 
is an $\alpha\in E_{k-1,k}$ and the other is a $\gamma\in E_{k+1,k+2}$, 
and the path $\ovl{\alpha\beta\gamma}$ is $f$-admissible.
	\item[{(b)}] There exists an integer $\chi(f)\in\5{0,1,2,3}$ 
such that for all values of $k$, the set $E_{k,k+1}$ contains exactly 
$\chi(f)$ yellow edges and the same number of red edges. 
	\end{Corollary}

\Pf (a): Consider an edge \6{\beta}{E^+_{k,k+1}}. By Lemma~\ref{L2}(b), 
$\beta$ cannot be coincident to another element of $E^+_{k,k+1}$.  By 
Lemma~\ref{L2}(d), there is a unique edge \6{\alpha}{E^+_{k-1,k}} meeting 
$\beta$.  By Lemma~\ref{L2}(c)(b), there is a unique edge 
\6{\gamma}{E^+_{k+1,k+2}} adjacent to $\beta$.  (Again, the case of $E^-$ 
is symmetric.)

(b): By (a), the set $E^+_{k,k+1}$ has at most $3$ elements (as including 
no coincidences) and there is a one-to-one correspondence between the 
elements of $E^+_{k,k+1}$ and $E^+_{k+1,k+2}$ for every $k$ (namely,  
$\alpha\leftrightarrow\beta\leftrightarrow\gamma$)  Now it is enough to 
use Lemma~\ref{L2}(a) and the fact that the cyclic graph $C_N$~is 
connected.~\fP

\begin{Definition} The number 
	\[ \chi(f) = |E^+_{k,k+1}| = |E^-_{k,k+1}| 
	\]
(as in Corollary~\ref{C1}(b)) will be called the 
\underline{\em characteristic number\/} of the winning strategy $f$.  $\Box$
\end{Definition}

The case $\chi(f)=1$ can be excluded outright, since it would imply the 
existence of an $f$-admissible $N$-periodic paths of both directions (red 
and yellow).  Now only three cases remain: $\chi(f)=0,2,$ and $3$.


\subsection{The case \mbox{$\chi(f)=0$}}\label{Sec0}  
Here we additionally suppose that $N\geq 5$ (the cases of $N=3,4$ being 
already settled).

In the case of $\chi(f)=0$ all the edges of $G$ are blue.  That was  
possible for $N=3$ and $N=4$ as shown in Section~\ref{Sec34}. But  
supposing $N\geq 5$ we are going to prove that $f$ would in fact  
be a losing strategy, implying $\chi(f)\neq 0$ for $N>4$. 

Take any edge $\ovl{ab}$ of graph $G_N$.  It has, in particular, $32$ 
different $f$-admissible extensions of length $N+1$, by $N-3$ edges to the 
left and $3$ edges to the right, of the form
	\begin{equation} \label{chi0}
\ovl{a_{ij}\,a_j\,u_1\,\ldots\,u_{N-4}\,a\,b\,b_p\,b_{pq}\,b_{pqr}} 
	\end{equation}
for $i,j,p,q,r\in\5{1,2}$, where the choice of vertices 
$u_1\ldots,u_{N-4}$ is fixed.  Observe that, while there are exactly $4$ 
edges of the form $\ovl{a_{ij}a_j}$ and $\ovl{b_pb_{pq}}$ alike, there may 
be either $2$ or $3$ vertices $b_{pq}$ (and $a_{ij}$ alike).

First, suppose vertex $b_{pq}$ assumes three different values. Then there 
are at least six edges $\ovl{b_{pq}b_{pqr}}$, while the number of the 
edges $\ovl{a_{ij}a_j}$ is four, making the path close as $6+4>9$.  Next, 
suppose there are just two different vertices $b_{pq}$.  We can make index 
$q$ point to these vertices, so that $b_{0q}=b_{1q}$ for $q=0,1$.  Since 
the $a_{ij}$ and $b_{pq}$ are both in a $3$-element set $V_m$, one can now 
fix $i,j,q$ so that $a_{ij}=b_{0q}=b_{1q}$.  Then, one of two paths 
$\ovl{b_pb_{pq}a_j}$ \9{p=0,1} must be admissible since 
$\ell_-(\ovl{b_{pq}a_j})=2$ and $b_p$ takes $2$ values.  Thus, the path 
(\ref{chi0}) acquires a closure with no use of vertex $b_{pqr}$.  (But in 
fact, $a_j=b_{pqr}$ for some $r$.)

	\begin{Corollary}\label{C0}
For $N>4$, every winning strategy $f$ has $\chi(f)\neq 0$.  
	\end{Corollary}


\subsection{The case \mbox{$\chi(f)=3$}}  \label{Sec3}
By Corollary~\ref{C1}, the 
yellow edges of graph $G$ are arranged as follows: 
	\[
	\begin{array}{l}
u_1^0 \xto u_1^1 \xto u_1^2 \xto\cdots\xto u_1^{N\!-\!1} \xto 
u_{\sigma(1)}^0 \xto u_{\sigma(1)}^1 \xto u_{\sigma(1)}^2 
	\cdots \\  
u_2^0 \xto u_2^1 \xto u_2^2 \xto\cdots\xto u_2^{N\!-\!1} \xto 
u_{\sigma(2)}^0 \xto u_{\sigma(2)}^1 \xto u_{\sigma(2)}^2 
	\cdots \\ 
u_3^0 \xto u_3^1 \xto u_3^2 \xto\cdots\xto u_3^{N\!-\!1} \xto 
u_{\sigma(3)}^0 \xto u_{\sigma(3)}^1 \xto u_{\sigma(3)}^2  
	\cdots  
	\end{array}
	\]
where $\5{u_1(k),u_2(k),u_3(k)} = V_k$ for all $k$ and 
\7{\sigma}{\5{1,2,3}}{\5{1,2,3}} is a permutation. 

(If $\sigma$ had a fixed point, then a yellow cycle of period~$N$ would 
defeat strategy~$f$.  Hence, $\sigma$ must be a rotation:
$\sigma(i)\equiv i\pm 1\,(\r{mod}\,3)$ for  $i=1,2,3$.  This observation 
will be used later to construct winning strategies.) 

Now let us locate the other colours.  By Corollary~\ref{C1}, the set 
$E_{1,2}$ contains three disjoint red edges.  Thus, we may assume
	\[ 
E^-_{1,2} = 
\5{\ovl{u_1(1)u_3(2)},\ \ovl{u_2(1)u_1(2)},\ \ovl{u_3(1)u_2(2)}} 
	\] 
(the other possibility being symmetric: 
$u_1u_2,\,u_2u_3,\,u_3u_1$). Considering the set $E_{0,1}$ (to the left 
of $E_{1,2}$) we see that 
$\ovl{u_3(0)u_2(1)}$ is red.  Indeed, since 
$\ovl{u_3(0)u_3(1)}\in E^+$, we have 
$\ovl{u_3(0)u_3(1)u_1(2)}\not\in\AA(f)$, so 
$\ovl{u_3(0)u_2(1)u_1(2)}\in\AA(f)$, implying that $\ovl{u_3(0)u_2(1)}$ 
must be the left-directed left continuation of 
$\ovl{u_2(1)u_1(2)}\in E^-$ (by Lemma~\ref{L2}(c)).  It follows that 
the edges in $E^-_{0,1}$  
have the same arrangement as in $E^-_{1,2}$.  Similarly, the sets 
$E^-_{k,k+1}$ and $E^-_{k+1,k+2}$ have the same arrangement for all 
$k=0,\ldots,N-2$, so we may assume that 
	\[  E^-_{k,k+1} = 
	\5{
\ovl{u_1^k u_3^{k+1}},\, \ovl{u_2^k u_1^{k+1}},\, \ovl{u_3^k u_2^{k+1}}
	}  
\ \ \ \ (k=0,1,2,\ldots,N-1).
	\] 
All the remaining edges of $G$ must be blue (by Lemma~\ref{L2}(b)).  
Here, by Lemma~\ref{L2}(e), the periodic $f$-admissible paths are 
precisely the periodic ones of a fixed colour.  

Now consider all three one-colour paths of length $N$, starting at vertex 
$u_1(0)$ and going to the right (for the red one, this means going back).  
If $N$ is divisible by $3$, they all end at $u_{\sigma(1)}(0)$.  But if 
$N$ is not divisible by $3$, they end at three distinct vertices of $V_N$, 
one of which must be $v_1(0)$.  That makes one of the paths close to  
defeat the strategy, which is a contradiction.  
A significant part of the situation for $N=5$ is illustrated in the 
diagram.  

	\begin{figure} \label{f5}
	\begin{center} \setlength{\unitlength}{0.6ex} 
	\begin{picture}(102,22)(-1,-1)  \thinlines  
\multiput(0,19)(20,0){5}{\makebox(0,0)[b]{$1$}}  
\multiput(0, 9)(20,0){5}{\makebox(0,0)[b]{$2$}}  
\multiput(0,-1)(20,0){5}{\makebox(0,0)[b]{$3$}}  
     \put(100,19){\makebox(0,0)[b]{\boldmath $2$}}  
     \put(100, 9){\makebox(0,0)[b]{\boldmath $3$}}
     \put(100,-1){\makebox(0,0)[b]{\boldmath $1$}} 
	\multiput(2,20)(20,0){5}{\vector(1,0){16}}  
   \put(18, 2){\vector(-1,1){16}} 
   \put(38, 9){\vector(-2,-1){16}}   
   \put(58,19){\vector(-2,-1){16}}
   \put(78, 2){\vector(-1,1){16}}
   \put(98, 9){\vector(-2,-1){16}}
\thicklines 
\put( 2,19){\line(2,-1){16}}  \put(22, 9){\line(2,-1){16}}
\put(42, 2){\line(1,1){16}}
\put(62,19){\line(2,-1){16}}  \put(82, 9){\line(2,-1){16}}
	\end{picture} 
\caption{Part of a typical strategy of characteristic $3$ on $C_5$} 
	\end{center}
	\end{figure}
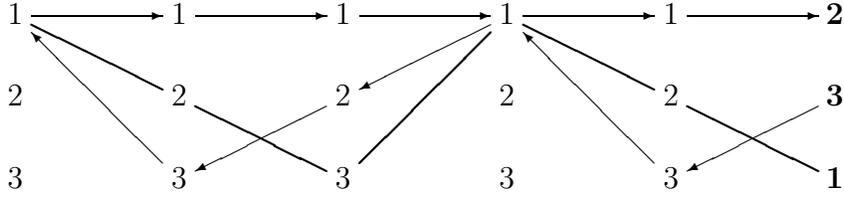


\subsubsection{Winning for \mbox{$3|N$}}    
If $N$ is a multiple of $3$ and strategy $f$ is colourable in the pattern 
just considered, then any one-colour path starting at $u_i(0)$ passes 
through $u_i(N)\neq u_i(0)$ and goes three times around the graph $G_N$, 
ending with period $3N\neq N$.  Moreover, the considered colour 
arrangement is always (for every $N$) given by some colourable strategy, 
defined in the following way.  Let
	\[
u_i(k+N) = u_{\sigma(i)}(k) ~~\r{for all}~~ k\in\Z,  
	\] 
where $\sigma$ is some fixed-point-free permutation, and 
	\begin{equation}\label{eq_chi3} 
f_k = \left[ \begin{array}{ccc} 2 & 1 & 1\\ 2 & 3 & 2\\ 3 & 3 & 1 
\end{array}\right]\ \ \ \ (k\in\Z), 
	\end{equation} 
using the convention: $f_k(i,j)$ in row~$i$, column~$j$.  
The definition is consistent since $f_k$ is rotation-invariant, i.e.,
	\[
f_k(\sigma(i),\sigma(j)) = \sigma(f_k(i,j)), 
	\] 
which can be checked directly.  (In fact, the strategies $f_k$ are 
uniquely determined by this rotation-invariant colouring, hence they must 
themselves be $\sigma$-invariant).

With this strategy, every directed edge is followed by an edge of the same 
direction, as in Lemma~\ref{L2}(c)).  Thus, any admissible periodic path 
has one colour, as in Lemma~\ref{L2}(e).  Since no admissible path has 
period $N$, the strategy is winning.


\subsubsection{The solution for \mbox{$\chi(f)=3$}}

	\begin{Corollary}\label{C3}
A winning strategy $f$ with $\chi(f)=3$ exists if and only if $N$ is 
divisible by~$3$.  If it exists, $f$ is unique up to isomorphism 
(induced by permutations of colours at the vertices).  
	\end{Corollary} 

{\em Remark.\/} The situation for $\chi(f)=3$ can be visualised on a 
torus obtained by rotating a triangle which at the same time makes $1/3$ 
of a full turn in its own plane.  This situation can also be viewed using 
a covering of graph $3\ast C_N$ by the graph $3\ast C_\infty$, where 
$C_\infty$ has edges between all pairs of consecutive integers and is the 
universal covering of $C_N$.  


\subsection{The case $\chi(f)=2$} Corollary~\ref{C1} implies the following 
arrangement of all the yellow edges and some blue edges of graph $G$:

\newlength{\utau}
\settowidth{\utau}{$u_{\tau(2)}^0$} 
\def\uuu#1{\makebox[\utau]{#1}}
	\[
	\begin{array}{l}
u_1^{0} \xto u_1^{1} \xto u_1^{2} \xto\cdots\xto u_1^{N\!-\!1} \xto 
u_{\tau(1)}^0 \xto u_{\tau(1)}^1 \xto u_{\tau(1)}^2 
	\cdots \\ 
u_2^0 \xto u_2^1 \xto u_2^2 \xto\cdots\xto u_2^{N\!-\!1} \xto 
u_{\tau(2)}^0 \xto u_{\tau(2)}^1 \xto u_{\tau(2)}^2 
	\cdots \\
u_3^0 \xline u_3^1 \xline u_3^2 \xline\cdots\xline u_3^{N\!-\!1} \xline  
\uuu{$u_3^0$} \xline \uuu{$u_3^1$} \xline \uuu{$u_3^2$} 
	\cdots  
	\end{array}
	\]
where $\5{u_1(k),u_2(k),u_3(k)} = V_k$ for all $k$ and 
\7{\tau}{\5{1,2}}{\5{1,2}} is a permutation. If $\tau$ were the identity, 
then two yellow cycles would have period $N$, contrary to the assumption 
that $f$ is winning.  Thus, $\tau$ must be a transposition: $\tau(1)=2$ 
and $\tau(2)=1$.  

Let us look for the red and the remaining blue edges.  By 
Lemma~\ref{L2}(b)(d), the red edges are crossing between 
rows of the yellow ones: 
	\[
E^-_{k,k+1} = \5{\ovl{u_1(k)u_2(k+1)},\, 
\ovl{u_2(k)u_1(k+1)}} ~~\r{for}~~  
k=0,1,\ldots,N-1 
	\]
and all the remaining edges are blue.  

Now if $N$ were an odd number, then the red (left-directed) path 
($f$-admissible by Lemma~\ref{L2}(e)) ending at vertex $u_1(0)$ would 
begin at vertex $u_{\tau(2)}(N)=u_1(0)$ and have period $N$, again 
contrary to the assumption that $f$ is winning.  Consequently, $N$ must be 
an even number.


\subsubsection{The strategy for \mbox{$\chi(f)=2$}} 

As in the case of $\chi(f)=3$, from the obtained colour arrangement one 
can deduce the form of strategy $f$.  One obtains:
	\begin{equation}\nonumber 
f_k = \left[ \begin{array}{ccc} 3 & 1 & 1\\ 2 & 3 & 2\\ 2 & 1  & 3 
\end{array}\right]\ \ \ \ (k\in\Z) 
	\end{equation} 
with the same convention as before, and the 
permutation $\sigma=(2,1,3)$.  Again, this colourable strategy $f$  
with single-colour admissiblility exists for all even numbers $N$ since 
$f_1$ is invariant under the permutation $(2,1,3)$. 

Since both directed paths have period $2N$, any $N$-periodic admissible 
path must be blue by Lemma~\ref{L2}(e).  The only admissible blue paths of 
length $2$ are:
	\[
\ovl{u_i^k u_3^{k+1} u_3^{k+2}},\  
\ovl{u_3^k u_3^{k+1} u_i^{k+2}},\  
\ovl{u_3^k u_i^{k+1} u_3^{k+2}}, ~~\r{where}~~ \6{i}{\5{1,2}}. 
	\]


\subsubsection{The question of winning with \mbox{$\chi(f)=2$}} 

Already for $N=2$ (despite the fact that $\chi(f)$ was not defined for 
$N<4$) a losing blue cycle can be observed, namely 
	\[
\ovl{\ldots v_3^0 v_1^1 v_3^2 v_2^3 v_3^4 v_1^5  \ldots}, 
	\]
where $v_1^1=v_2^3$. 
 
Now consider the case $N=4$.  Then any admissible path containing an 
edge $\ovl{v_3^k v_3^{k+1}}$ could not close with period $4$. At the same 
time, any admissible path containing no such edge must have the form 
$\ovl{\ldots31323132\ldots}$, i.e.~(up to a shift),
	\[
	\ldots 
\ovl{u_3(0)u_1(1)u_3(2)u_2(3)u_3(4)u_1(5)u_3(6)u_2(7)},  
	\ldots, 
	\]
which has period $8$ as $u_1(1)\neq u_1(5)$.  This shows that $f$ is a 
winning strategy for $N=4$.  

If, however, $2|N$ and $N\geq 6$, then 
there exists the following admissible blue path of period $N$:
	\[
\ldots \ovl{3313(31)(32)(31)(32)\ldots(3j)} \ldots 
	\] 
(where $j\in\5{1,2}$ and $j\equiv N/2\,(\r{mod}\,2)$).   
Consequently, the strategy $f$ is loosing for $N>4$.  The situation for 
$N=8$ is illustrated in the diagram.    

	\begin{figure} \label{f8}
	\begin{center} \setlength{\unitlength}{0.6ex}
	\begin{picture}(124,24)(-2,-2)  \thinlines 
   \multiput(0,19)(15,0){8}{\makebox(0,0)[b]{$1$}}  
   \multiput(0, 9)(15,0){8}{\makebox(0,0)[b]{$2$}}  
   \multiput(0,-1)(15,0){8}{\makebox(0,0)[b]{$3$}}  
\put(120,19){\makebox(0,0)[b]{\boldmath $2$}}
\put(120, 9){\makebox(0,0)[b]{\boldmath $1$}}
\put(120,-1){\makebox(0,0)[b]{\boldmath $3$}}
   \multiput(2,20)(15,0){8}{\vector(1,0){11}}  
   \multiput(2,10)(15,0){8}{\vector(1,0){11}}  
\multiput(13,18.67)(15,0){8}{\vector(-3,-2){11}} 
\multiput(13,11.33)(15,0){8}{\vector(-3, 2){11}}  
    \thicklines 
\put( 2,0){\line(1,0){11}} 	
	\put(17,1.33){\line(3, 2){11}} 
	\put(32,8.67){\line(3,-2){11}}    
\put(47,0){\line(1,0){11}} 
	\put(62,1.33){\line(3, 2){11}}    
	\put(77,8.67){\line(3,-2){11}}    
		\put( 92, 2.67){\line(3, 4){11}}    
		\put(107,17.33){\line(3,-4){11}}    
	\end{picture} 
\caption{A diagram of the typical strategy of characteristic $2$ on $C_8$, 
showing all the directed edges and a critical undirected path.} 
	\end{center}
	\end{figure}
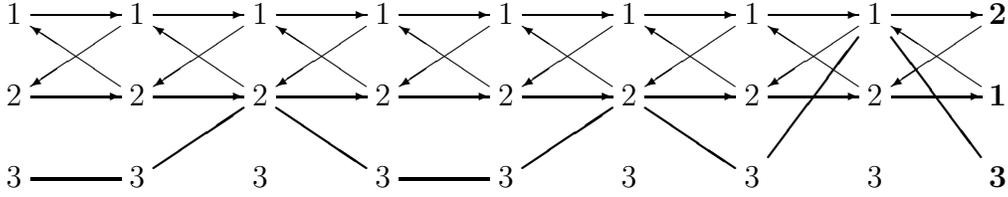


\subsubsection{The solution for \mbox{$\chi(f)=2$}}

	\begin{Corollary}\label{C2}
A winning strategy $f$ with $\chi(f)=2$ exists only for $N=4$ (and  
is unique up to isomorphism). 
	\end{Corollary}

{\em Remark.\/} The above situation for $2|N$ can be visualised as a 
M\"obious band with the yellow cycle on the boundary and the red cycle 
inside, completed with a separate blue cycle which is not admissible.  
The edges can be drawn on a Klein bottle arising from this construction.  
As before, this is equivalent to using an appropriate covering of graph 
$3\ast C_N$ by the graph $3\ast C_\infty$.

 
Thus we have proved Theorem~\ref{T1}.   \fP

\section{Corollaries}  

The configuration of three colours for each player is a maximal 
solvable one for the cycle graphs. 

\begin{Corollary}\label{C23}  The hat game on any cycle $C_N$ \9{n>4} with 
the height function $h$ satisfying $h(1)=4$ and $h(k)=3$ for $k=2,\ldots,N$
is losing, i.e., $\mu(h)=0$. 
\end{Corollary}

\Pf  If $f$ were a winning strategy for this game, then it would also be 
winning for any of its $3$-colour restrictions.  It follows that choosing 
any $k\in\5{1,2,3,4}$ and changing any values $f_1(i,j)=k$ into any values 
$f_1(i,j)\neq k$ would result in a winning strategy for the $3$-colour 
game.  By Corollary~\ref{C3}, such a strategy is unique up to permutations 
of colours.  Now, $f_1$ assumes some values, so for instance, we have 
$f_1(i,j)=4$ for some pair(s) \9{i,j}.  But the form (\ref{eq_chi3}) of 
the individual strategy shows that $f_1$ must assume each value three 
times.  Some arbitrary choice of $f_1(i,j)\neq 4$ could always change 
that, contrary to the fact that $f$ should remain a winning strategy.   
\fP 

Even if $N$ is not divisible by $3$, the probability of winning by 
using a {\em random\/} strategy equals $1-\9{\fr{2}{3}}^N$, but a much 
more effective non-random strategy can be chosen (assuming that the 
adversary does not know about it):

\begin{Corollary}\label{C_E} In the three-colour game on any cycle $C_N$ 
\9{3\not|N} there exists a strategy for which the probability of winning is 
$\geq 1-3^{-N+1}$.  \end{Corollary}

\Pf
This follows from the fact that the strategy described in Section~\ref{Sec3} 
has at most three admissible $N$-periodic paths. 
\fP  


\xxx{General graphs}

\paragraph{Acknowledgements.} My thanks go to Jaros\u.aw 
Grytczuk and S\u.awomir Por\e.bski for acquainting me with the hat 
problems, and Marcin Krzywkowski for reading the manuscript. 


\end{document}